\documentclass{article}
\usepackage{amsthm}
\usepackage{amssymb}
\usepackage{amsfonts}
\usepackage{color}
\usepackage{amsmath}
\begin{document}

\newcommand{\A}{\cal{A}}
\newcommand{\R}{\cal{R}}
\newcommand{\F}{\mathbb{F}}
\newcommand{\C}{\cal{C}}
\newcommand{\N}{\cal{N}}
\newcommand{\cH}{\cal{H}}
\newcommand{\G}{\cal{G}}
\newcommand{\Z}{\cal{Z}}
\newcommand{\T}{\cal{T}}
\newcommand{\Q}{\cal{Q}}
\newcommand{\cS}{\cal{S}}
\newcommand{\W}{\cal{W}}
\newcommand{\E}{\cal{E}}
\newcommand{\D}{\cal{D}}
\newcommand{\cP}{\cal{P}}
\newcommand{\cL}{\cal{L}}
\newcommand{\M}{\cal{M}}
\newcommand{\V}{\cal{V}}
\newcommand{\B}{\cal{B}}
\newcommand{\U}{\cal{U}}
\newcommand{\K}{\cal{K}}
\newcommand{\CC}{{\mathbb C}}
\newcommand{\Aa}{{\mathbb A}}
\newcommand{\FF}{{\mathbb F}}
\newcommand{\PP}{{\mathbb P}}
\newcommand{\ov}{\overline}
\newcommand{\eps}{\varepsilon}
\newcommand{\wh}[1]{\widehat{#1}}
\newtheorem{theorem}{Theorem}[section]
\newtheorem{proposition}[theorem]{Proposition}
\newtheorem{lemma}[theorem]{Lemma}
\newtheorem{corollary}[theorem]{Corollary}
\newtheorem{remark}[theorem]{Remark}
\newtheorem{definition}[theorem]{Definition}            
\newtheorem{example}[theorem]{Example}            
\newtheorem{conjecture}[theorem]{Conjqecture}
\newtheorem{question}[theorem]{Question}

\title{The study of symmetries: some general techniques\thanks{The author was supported by grants P1-0288 and N1-0368 from ARIS, Slovenia.}}
\author{Peter \v Semrl\footnote{Institute of Mathematics, Physics, and Mechanics, Jadranska 19, SI-1000 Ljubljana, Slovenia, peter.semrl@fmf.uni-lj.si; Faculty of Mathematics and Physics, University of Ljubljana, Jadranska 19, SI-1000 Ljubljana, Slovenia}
        }

\date{}
\maketitle

\begin{abstract}
Let $S(H)$ be the set of all self-adjoint bonded linear operators on $H$ and $\mathcal{V} \subset S(H)$ a subset that is pertinent in mathematical foundations of quantum mechanics. A symmetry is a bijective map $\phi :\mathcal{V} \to \mathcal{V}$ which is an automorphism with respect to one or more relations and/or operations on $\mathcal{V}$ that are relevant in mathematical physics. We will explain several ideas that can be used when studying the general form of symmetries.
\end{abstract}
\maketitle



\section{Introduction}

Let $H$ be a complex Hilbert space, $\dim H \ge 2$. We denote by $B(H)$ and $S(H)\subset B(H)$ the algebra of all bounded linear operators on $H$ and the subset of all self-adjoint operators, respectively. In mathematical foundations of quantum mechanics bounded observables are represented by self-adjoint operators. There are several subsets of $S(H)$ that are important in mathematical formalization of quantum mechanics. 
Let $A \in S(H)$. Then $A \ge 0$ if $\langle Ax , x \rangle \ge 0$ for every $x \in H$. Such operators are called positive operators. When $A \ge 0$ and $A$ is invertible we say that $A$ is a positive definite operator and write $A >0$. The set
of all positive operators $S_{\ge} (H) = \{ A \in S(H)\, : \, A \ge 0\}$ and the set of all positive definite operators $S_> (H) = \{ A \in S(H)\, : \, A > 0\}$ play  important roles in physics. The effect algebra $E(H) = \{ A \in S(H)\, : \, A \ge 0 \ \, {\rm and} \ \,  I-A\ge 0 \}$
is of basic importance  in Ludwig's axiomatic formulation of quantum mechanics, see \cite{BLPY}, \cite{Dav}, \cite{Kraus}, \cite{LudI}, and \cite{LudII}. We denote by $P(H)$ the set of all projections on $H$, $P(H) = \{ A \in S(H)\, : \, A^2 = A = A^\ast \}$. Pure states are repreented by projections of rank one. We denote by $P_1 (H) = \{ P \in P(H)\, : \, \dim {\rm Im}\,( P) = 1 \}$ the set of all pure states. 

On the set $S(H)$ or its subsets we can define several operations and relations that are relevant in mathematical foundations of quantum mechanics. Let $A,B \in S(H)$. Then $A \le B$ if and only if $B-A \ge 0$, and similarly, $A < B \iff B-A > 0$.
For bounded observables $A,B \in S(H)$ we have $A \le B$ if and only if the mean value (expectation) of $A$ in every state is less than or equal to the mean value of $B$ in the same state. The relation $\le$ on $S(H)$ is often called Loewner order.
Two bounded observables are compatible if they can be measured jointly. If $A,B$ are the corresponding self-adjoint operators then the two bounded observables are compatible if and only the operators $A$ and $B$ commute. The situation is more involved in the case of effect algebras.
Two effects that can be measured together by applying a suitable apparatus are called coexistent effects. When we switch to the
language of mathematics this translates into the following
definition: Two effects $A,B \in E(H)$ are said to be coexistent if there exist effects $E,F,G \in E(H)$ such that
$$
A = E + G  \ \ \ {\rm and} \ \ \ B = F + G \ \ \ {\rm and} \ \ \ E + F + G \in E(H).
$$
For any effect $A$ we denote by $A^\perp$ the orthocomplent of $A$, $A^\perp = I -A$.
For two self-adjoint operators $A$ and $B$ the product $AB$ is in not always self-adjoint. But $S(H)$ is closed for the Jordan product $A \circ B = (1/2)(AB + BA)$ as well as for the triple Jordan product $A \circ B = ABA$. In the setting of effect algebras a more natural concept is the sequential product defined by $A \circ B = A^{1/2} B A^{1/2}$. On the set of positive operators an important operation is the geometrical mean defined by $A \sharp B = A^{1/2} (A^{-1/2} B A^{-1/2} )^{1/2} A^{1/2}$. Transition probability between pure states $P,Q \in P_1 (H)$ is defined as ${\rm tr}\, (PQ)$, where ${\rm tr}\, A$ stands for the trace of matrix $A$. In particular, the transistion probability between $P$ and $Q$ is zero if
 ${\rm tr}\, (PQ) = 0$, which is easily seen to be equivalent to $PQ = QP = 0$. In this case we write $P \perp Q$. It is also easy to see that $P \perp Q$ if and only if $P \not= Q$ and $P$ and $Q$ commute. 
We can extend the notion of orthogonality to arbitrary projections. For $P,Q \in P(H)$ we write $P \perp Q$ if the images of $P$ and $Q$ are orthogonal which is equivalent to $PQ=QP=0$.

There are many more subsets of $S(H)$ and relations and operations on these subsets that are relevant in mathematical physics. Let $\mathcal{V}$ be a subset of $S(H)$ and $\circ$ a binary operation on $\mathcal{V}$. We say that a map 
$\phi : \mathcal{V} \to \mathcal{V}$ preserves the operation $\circ$ if
$$
\phi (A \circ B) = \phi (A) \circ \phi (B), \ \ \ A,B \in \mathcal{V}.
$$
Assume that $\sim$ is a relation  on $\mathcal{V}$. Then a map 
$\phi : \mathcal{V} \to \mathcal{V}$ preserves the relation $\sim$ if for every pair $A,B \in \mathcal{V}$ we have
$$
A \sim B \iff \phi (A) \sim \phi (B) .
$$
A bijective map $\phi : \mathcal{V} \to \mathcal{V}$ that preserves a certain operation $\circ$ is called a symmetry of $\mathcal{V}$ with respect to the operation $\circ$. Similarly, a bijective map $\phi : \mathcal{V} \to \mathcal{V}$ preserving $\sim$ is called a symmetry of $\mathcal{V}$ with respect to the relation $\sim$. Sometimes we are interested in symmetries that preserve more than just one operation or relation. For example, a bijective map $\phi : E(H) \to E(H)$ with the property that for every pair $A,B \in E(H)$ we have
\begin{equation}\label{less}
A \le B \iff \phi (A) \le \phi (B)
\end{equation}
and
\begin{equation}\label{ortho}
\phi (A^\perp) = \phi (A)^\perp
\end{equation}
is called a symmetry of the effect algebra with respect to Loewner order and orthocomplementation. We could have defined symmetries shortly: We assume that $\mathcal{V}$ is a subset of $S(H)$ and that on this subset we have certain relations and/or certain operations. A symmetry of $\mathcal{V}$ with respect to these relations and operations is an automorphism of $\mathcal{V}$ equipped with these relations and operations. In particular, a bijective map $\phi : E(H) \to E(H)$ satisfying (\ref{less}) and (\ref{ortho}) is called an ortho-order automorphism of the effect algebra $E(H)$.

Let us give a few examples of symmetries. First we need to recall a definiton of a conjugate-linear operator on $H$. A map $T : H \to H$ is a bounded conjugate-linear operator if it is continuous and for every pair $x,y \in H$ and every complex number $\lambda$ we have $T(x+y) = Tx + Ty$ and $T(\lambda x) = \overline{\lambda} Tx$. For every bounded conjugate-linear operator $T : H \to H$ there esists a unique bounded conjugate-linear operator $T^\ast : H \to H$ such that for every pair $x,y \in H$ we have
$$
\langle Tx,y \rangle = \overline{ \langle x, T^\ast y \rangle} = \langle  T^\ast y, x \rangle .
$$
In the case that a conjugate-linear operator $T$ is a bijective isometry we say that it is an antiunitary operator, that is, a conjugate-linear operator $U : H \to H$ is antiunitary if it is bijective and 
$$
\langle Ux, Uy \rangle = \langle y, x \rangle, \ \ \ x,y \in H.
$$

Let $T : H \to H$ be an invertible bounded linear or conjugate-linear operator. A map $\phi : S_{\ge} (H) \to S_{\ge} (H)$ defined by
\begin{equation}\label{pos}
\phi (A) = TAT^\ast , \ \ \ A \in S_{\ge} (H),
\end{equation}
is a symmetry of $S_{\ge} (H)$ with respect to Loewner order. Similarly, the congruence transformation  $\phi : S_{>} (H) \to S_{>} (H)$ defined by
\begin{equation}\label{posdef}
\phi (A) = TAT^\ast , \ \ \ A \in S_{>} (H),
\end{equation}
is an order automorphism of $S_{>} (H)$. Let further $S \in S(H)$ be any self-adjoint operator. Then the map $\phi : S(H) \to S(H)$ given by
\begin{equation}\label{vse}
\phi (A) = TAT^\ast  + S, \ \ \ A \in S (H),
\end{equation}
is a symmetry of $S(H)$ with respect to Loewner order. Let $U : H \to H$ be a unitary or antiunitary operator. Then the unitary-antiunitary similarity $\phi : E(H) \to E(H)$
defined by 
\begin{equation}\label{mija}
\phi (A) = UAU^\ast, \ \ \ A \in E(H),
\end{equation}
is an order automorphism of the effect algebra.

Similarly, every unitary-antiunitary similarity $A \mapsto UAU^\ast$, $A \in S(H)$, is a symmetry of $S(H)$ with respect to commutativity (=compatibility in the language of quantum mechanics) and also with respect to Jordan product. Every unitary-antiunitary similarity is a symmetry of $P_1(H)$ with respect to the transition probability, that is, for every unitary or antiunitary operator $U : H \to H$ the map $\phi : P_1(H) \to P_1 (H)$ given by $\phi (P) = UPU^\ast$, $P \in P_1 (H)$, is bijective and
\begin{equation}\label{mcda}
{\rm tr}\, (\phi (P) \phi (Q)) = {\rm tr}\, (PQ), \ \ \ P,Q \in P_1 (H).
\end{equation}

The usual goal when studying symmetries is to describe the general form of these maps. Above we have given a few easy examples of symmetries. It is somewhat surprising that in many cases there are no other symmetries. For example, if $\phi$ is a symmetry of $S_{\ge} (H)$ with respect to Loewner order then it has to be of the form (\ref{pos}) for some invertible bounded linear or conjugate-linear operator $T : H \to H$, and similarly,  if $\phi$ is a bijective map of $S_{>} (H)$ onto itself which preserves Loewner order then there exists an invertible bounded linear or conjugate-linear operator $T : H \to H$ such that (\ref{posdef}) holds \cite{Mol}. Moreover, every symmetry of $S(H)$ with respect to $\ge$ is of the form (\ref{vse}) \cite{Mo1}. Every ortho-order automorphism of the effect algebra $E(H)$ is a  unitary-antiunitary similarity \cite{LudI}.  And finally, if $\phi : P_1(H) \to P_1(H)$ is a bijective map such that for every pair $P,Q \in P_1(H)$ we have (\ref{mcda}) then there exists a unitary or antiunitary operator $U : H \to H$ such that $\phi P) = UPU^\ast$, $P \in P_1 (H)$. This is the famous Wigner's unitary-antiunitary theorem, one of the most important results in the mathematical foundations of quantum mechanics.

In some other cases the structure of symmetries is much more complicated. For example, if $A,B \in S(H)$ commute then $A^3 - 5A + 2I$ and $B^2$ commute as well. This quickly implies that besides unitary-antiunitary similarities there are other symmetries of $S(H)$ with respect to commutativity. At this point the natural conjecture is the following: Let $\phi : S(H) \to S(H)$ be a bijective map with the property that for each pair of operators $A,B \in S(H)$ we have
$$
AB = BA \iff \phi (A) \phi (B) = \phi (B) \phi (A).
$$
Then there exists a unitary or antiunitary operator $U : H \to H$ and for every $A \in S(H)$ there is a real valued bounded Borel function $f_A$ defined on the spectrum of $A$ such that
$$
\phi (A) = U f_A (A) U^\ast .
$$
And indeed, this conjecture was confirmed in the case when $H$ is a separable Hilbert space with $\dim H \ge 3$, see \cite{MoS}.

The following example given in
\cite{MoK} shows that the structure of symmetries of $E(H)$ with respect to Loewner order is quite complicated. 
Let $T \in E(H)$ be invertible. Then the map
\begin{equation}\label{vau}
A \mapsto 
\left( { T^2 \over 2I - T^2 }\right)^{-1/2} \left( (I-T^2 + T(I+A)^{-1}T)^{-1} - I \right) 
\left( { T^2 \over 2I - T^2 }\right)^{-1/2}
\end{equation}
is an automorphism of $(E(H), \le )$. In order to verify this
we recall the well-known fact that if $A,B > 0$, then
we have $A\le B$ if and only if $B^{-1} \le A^{-1}$. 
Further, for $E,F \in S(H)$ with $E < F$ we define 
 $[E,F] = \{ C \in S(H)\, : \, E \le C \le F\}$. Let $E_j, F_j \in S(H)$ satisfy $E_j < F_j$, $j=1,2$. A  bijective map $\phi : [E_1,F_1] \to [E_2 , F_2]$ is called
an order anti-isomorphism if for every $A,B \in [E_1,F_1]$ we have $A \le B \iff \phi (B) \le \phi (A)$.

Let $S\in S(H)$ be any self-adjoint operator. Then the transformation $$A \mapsto A+S$$ is an order isomorphism of
$[E,F]$ onto $[E+S, F+S]$. Further, if $S: H \to H$ is any bounded invertible linear or conjugate-linear map,
then the map $$A \mapsto SAS^*$$ is an order isomorphism of $[E,F]$ onto $[SES^* , SFS^* ]$. And finally, assume that $0 < E < F$. Then
 the transformation $$A\mapsto A^{-1}$$
is an order anti-isomorphism of $[E,F]$ onto $[F^{-1}, E^{-1}]$.

It is now not difficult to understand the above example. The map $\tau$ given by
$$
A \mapsto 
\tau (A) = \left( I-T^2 + T(I+A)^{-1}T \right)^{-1} - I  
$$
is easily seen to be a compositum of a few order isomorphisms and exactly two order anti-isomorphisms as follows:
$$
A \mapsto I + A \mapsto (I + A )^{-1} \mapsto T (I+A)^{-1}T \mapsto
$$
$$
I - T^2 + T (I+A)^{-1}T \mapsto (I - T^2 + T (I+A)^{-1}T)^{-1} \mapsto 
$$
$$
(I - T^2 + T (I+A)^{-1}T)^{-1} - I.
$$
Hence, $\tau$ is an order isomorphism  of $E(H)$ onto $[\tau (0), \tau (I)]$.
Obviously, $\tau (0) = 0$ and $\tau (I) = T^2 (2I -T^2 )^{-1}$. The map (\ref{vau}) can be rewritten as
$$
A \mapsto 
\left( { T^2 \over 2I - T^2 }\right)^{-1/2} \, \tau (A) \,
\left( { T^2 \over 2I - T^2 }\right)^{-1/2}
$$
and is therefore an order automorphism of $E(H)$. Let us just mention that the problem of describing the general form of symmetries of effect algebras with respect to Loewner order has been solved in \cite{MS} and \cite{Se0}.

By now a lot of results have been proved on the structure of symmetries. At first glance the proofs of different results seem to be quite unrelated. Still, there are some general ideas that were used repeatedly in this active research area. It is the aim of this expository paper to explain three approaches that might be used when trying to solve a problem of desribing the general form of symmetries. It should be mentioned that the techniques we will explain are quite often just the first step in the study of certain symmetries. Thus, the ideas that we will present might be of  some help when one starts to work on such a problem but in each specific case quite a lot of additional work might be needed to solve the problem completely.

\section{A method based on adjacency preservers}

Two operators $A,B \in S(H)$ are said to be adjacent if the difference $B-A$ is an operator of rank one.
Quite a few results describing the general form of symmetries were proved by reducing the original problem to the problem of describing the general form of maps preserving adjacency, see \cite{Se2}, \cite{Se1}, and the references therein.
We want to present here just the main idea without paying too much attention to technical details. In order to make our explanation as simple as possible we restrict ourselves to the finite-dimensional case. An interested reader can see what kind of ideas are needed to extend our method from the finite-dimensional case to the infinite-dimensional setting in \cite{Se2}. 

When we are dealing with self-adjoint operators on a Hilbert space $H$ with $\dim H = n < \infty$, we identify them with $n \times n$ hermitian matrices. We denote the set of all $n \times n$ hermitian matrices by $H_n$.
We need the following result from \cite{HuS}.

\begin{theorem}\label{ftghm}
Let $n\ge 2$ be an integer and $\phi : H_n \to H_n$ a map such that
for every pair $A,B \in H_n$ the matrices $A$ and $B$ are adjacent if and only if
$\phi (A)$ and $\phi (B)$ are adjacent. Then there exist $c\in \{ -1, 1 \}$,
an invertible $n\times n$ complex matrix $T$, and $S\in H_n$ such that either
$$
\phi (A) = cTAT^* + S, \ \ \ A\in H_n,
$$
or
$$
\phi (A) = cTA^t T^* + S, \ \ \ A\in H_n.
$$
Here, $A^t$ denotes the transpose of $A$.
\end{theorem}

Actually, the main theorem in \cite{HuS} is much stronger but the above weaker version is all that we need. It should be mentioned that the study of adjacency preserving maps was initiated by L. K. Hua in a series of papers \cite{Hu1} - \cite{Hu8}.

We will start by sketching the proof of the structural theorem for symmetries of $H_n$ with respect to Loewner order $\le$.

\begin{theorem}\label{orderpres}
Let $n\ge 2$ be an integer and $\phi : H_n \to H_n$ a bijective map such that
for every pair $A,B \in H_n$ we have
$$
A \le B \iff \phi (A) \le \phi (B).
$$
Then there exist 
an invertible $n\times n$ complex matrix $T$ and $S\in H_n$ such that either
$$
\phi (A) = TAT^* + S, \ \ \ A\in H_n,
$$
or
$$
\phi (A) = TA^t T^* + S, \ \ \ A\in H_n.
$$
\end{theorem}

The idea that we  would like to explain in this section is to reduce a given problem concerning symmetries to the problem of describing maps preserving adjacency. Once we do such a reduction we can apply Theorem \ref{ftghm}.
In our special case we will characterise the adjacency relation by the relation $\le$. So, if our map preservers $\le$ then it has to preserve adjacency and therefore it has to be of the nice form described in Theorem \ref{ftghm}. 

After this short explanation we start with the sketch of the proof of Theorem \ref{orderpres}. Assume first that $A,B \in H_n$ are adajcent. Then $B-A$ is of rank one, that is, there exist a rank one projection $P \in H_n$ and a nonzero real number $t$ such that $B-A = tP$. Clearly, if $t >0$, then $B \ge A$, and if $t$ is negative, then $B \le A$. Thus, $A$ and $B$ are comparable, that is, $A \le B$ or $B \le A$. Let us consider only the case that $t >0$. Take any pair $C,D \in H_n$ satisfying $A \le C,D \le B = A +tP$. It is easy to see that $C = A + pP$ and $D = A + qP$ for some real numbers $p,q \in [0,t]$. Clearly, we have $C \le D$ or $D \le C$ depending on which of the real numbers $p,q$ is larger. We have started with the assumption that $A$ and $B$ are adjacent and we have seen that in this case the following is true:
\begin{itemize}
\item[$( \dagger )$] $A$ and $B$ are comparable and any two hermitian matrices $C,D$ that are in between $A$ and $B$ are comparable.
\end{itemize}
Let us now consider $A,B \in H_n$, $A \not= B$, such that $A$ and $B$ are not adjacent. We will show that the condition $( \dagger )$ is not satisfied. If $A$ and $B$ are not comparable we are done. So, we may assume that $A \le B$ or $B \le A$. Again we will consider just the first possibility. Thus, $B-A$ is positive and is of rank at least two. It follows that there exist rank one projections $P,Q$, $P\not= Q$, and positive real numbers $p,q$ such that $pP, qQ \le B-A$. Note that neither $pP \le qQ$ nor $qQ \le pP$. Clearly, $A \le A + pP, A + qQ \le B$ and $A + pP, A + qQ$ are not comparable, and therefore $A$ and $B$ do not satisfy $( \dagger )$.

We have shown that if $A,B \in H_n$ and $A\not= B$, than $A$ and $B$ are adjacent if and only if $(\dagger)$ holds true. Let $\phi : H_n \to H_n$ be a bijective map such that
for every pair $A,B \in H_n$ we have
$A \le B \iff \phi (A) \le \phi (B)$. Then it is clear that for every pair of matrices $A,B \in H_n$, $A \not=B$, the pair $A,B$ satisfies $(\dagger)$ if and only if the pair $\phi (A), \phi (B)$ satisfies $(\dagger)$. This implies that 
for every pair $A,B \in H_n$ the matrices $A$ and $B$ are adjacent if and only if
$\phi (A)$ and $\phi (B)$ are adjacent. By Theorem \ref{ftghm} there exist $c\in \{ -1, 1 \}$,
an invertible $n\times n$ complex matrix $T$, and $S\in H_n$ such that either
$\phi (A) = cTAT^* + S$ for every $A\in H_n$,
or $\phi (A) = cTA^t T^* + S$ for every  $A\in H_n$. Applying the fact that $\phi$ preserves Loewner order we see that in both cases $c=1$. This completes the proof.

A similar idea works when we want to characterize symmetries of $H_n$ with respect to the triple Jordan product. 
Our goal is to verify the following statement.

\begin{theorem}\label{triple1}
Let $\phi : H_n \to H_n$ be a bijective 
map with the property that for every pair $A,B \in H_n$ we have
$$
\phi (ABA) = \phi (A) \phi (B) \phi (A). 
$$
Then there exist $c \in \{ -1, 1 \}$ and a unitary $n \times n$ matrix $U$ such that either
$$
\phi (A) = cUAU^*
$$
for every $A\in H_n$, or
$$
\phi (A) = cUA^t U^*
$$
for every $A\in H_n$.
\end{theorem}

The main steps of the proof are as follows. It is not difficult to check that $\phi (0) =0$.
It is also rather eays to check that for a nonzero $A \in H_n$ the following are equivalent:
\begin{itemize}
\item $A$ is a rank one matrix,
\item if $B \in H_n$ and
$\{ BCB \, : \, C \in H_n \}$ is a proper subset of $\{ ACA \, : \, C \in H_n \}$, then $B=0$.
\end{itemize}
It follows that $\phi$ maps the set $H_{n}^1$ of rank one matrices onto itself.  For $A,B \in H_{n}^1$ we write $A \perp B$ if the images of $A$ and $B$ are orthogonal one-dimensional spaces. This is equivalent to $ABA = 0$.
Thus, for every pair of rank one matrices
$A,B$ we have $A\perp B$ if and only if $\phi (A) \perp \phi (B)$.

Let $y$ be a nonzero vector.
We denote by ${\cal R}_y \subset H_{n}^1$ the set of all rank one matrices whose images are orthogonal to $y$.
It is clear that for every nonzero vector $y$ we
can find a nonzero vector $z$ such that $\phi ( {\cal R}_y ) = {\cal R}_z$. It is not difficult to verify that for an arbitrary pair $A,B \in H_n$ the matrices $A$ and $B$ are adjacent if and only if 
there exists a nonzero vector $y$ such that $\{ C\in H_{n}^1  \, : \, CAC = CBC \} = {\cal R}_y$. 
Once we have this charaterization of adjacency involving Jordan triple product we immediately conclude that 
$\phi$ preserves adjacency and then it is easy to complete the proof by applying Theorem \ref{ftghm}.

Much more examples illustrating the efficency of the method illustrated by the above two examples can be found in \cite{Se2, Se1}.

\section{Applying projective geometry to describe the structure of symmetries}

Let $H$ be a Hilbert space. Throughout this section we will assume that $\dim H \ge 3$. For any nonzero $x \in H$ we denote by $[x]$ the one-dimensional subspace spanned by $x$. The projective space $\PP (H)$ is the set of all one-dimensional subspaces of $H$, $\PP (H) = \{ [x] \, : \, x \in H \setminus \{ 0 \} \}$. Let $\xi : \PP (H) \to \PP (H)$ be a bijective map such that for every triple  $x,y,z \in H \setminus \{ 0 \}$ we have
\begin{equation}\label{kajmj}
 [x] \subset  [y] +  [z] \iff   \xi ([x] )\subset \xi( [y] ) + \xi ( [z] ).
\end{equation}
The fundamental theorem of projective geometry tells that then there exists a bijective semilinear map $T : H \to H$ such that $\xi ([x] ) = [Tx ]$, $ x \in H \setminus \{ 0 \}$, see \cite{Fa} for a more general result with a simple proof. Recall that a map $T: H \to H$ is said to be semilinear if there exists a field automorphism $f: \mathbb{C} \to \mathbb{C}$ such that for all $x,y \in H$ and $\lambda \in \mathbb{C}$ we have
$$
T(x+y) = Tx + Ty \ \ \ {\rm and} \ \ \ T(\lambda X) = f(\lambda) Tx.
$$
The special cases of semilinear maps are linear maps and conjugate-linear maps. But there exist automorphisms of the complex field different from the identity and the complex conjugation and therefore there exist semilinear maps that are neither linear nor conjugate-linear. Clearly, we can identify $\PP (H)$ with $P_1(H)$ in a natural way, that is, we identify   $[x]$, $x\not=0$, with the rank one projection onto $[x]$.

As already mentioned the famous Wigner's theorem states that for every bijective map $\phi : P_1(H) \to P_1 (H)$ satisfying (\ref{mcda}) there exists a unitary or antiunitary operator $U : H \to H$ such that $\phi (P) = UPU^\ast$
for every $P \in P_1 (H)$. Uhlhorn's improvement of Wigner's unitary-antiunitary theorem states that one gets the same conclusion under the weaker assumption that only the zero transition probability is preserved. More precisely, we have the following theorem.

\begin{theorem}\label{uhlhorn}
Let $H$ be a Hilbert space with $\dim H \ge 3$ and $\phi : P_1(H) \to P_1 (H)$ a bijective map such that for every pair $P,Q \in P_1 (H)$ we have
\begin{equation}\label{hrvat}
PQ=0 \iff \phi (P) \phi (Q) = 0.
\end{equation}
Then 
there exists a unitary or antiunitary operator $U : H \to H$ such that $\phi (P) = UPU^\ast$
for every $P \in P_1 (H)$.
\end{theorem}

We will see that this theorem follows quite trivially from the fundamental theorem of projective geometry. For any subset $\mathcal{P} \subset P_1 (H)$ we denote by $\mathcal{P}^\perp$ the orthogonal complement of $\mathcal{P}$ in $P_1(H)$, that is,
$$
\mathcal{P}^\perp = \{ Q \in P_1(H) \, : \, PQ = 0 \ \, {\rm for} \ \, {\rm every} \ \, P \in \mathcal{P}  \}.
$$
Clearly, for every subset $\mathcal{P} \subset P_1 (H)$ we have $\phi ( \mathcal{P}^\perp ) = \phi (\mathcal{P})^\perp$.
When identifying $P_1(H)$ with $\PP (H)$ as above we see that $\phi$ induces in a natural way a map $\xi : \PP (H) \to \PP (H)$. For $x,y \in H$, $x,y \not= 0$, we have $\xi ( [x] ) =  [y]$ if and only if $\phi (P) = Q$, where $P$ and $Q$ are rank one projections onto $ [x]$ and $ [y]$, respectively. Let $x,y,z \in H$ be nonzero vectors and $P,Q,R$ the rank one projections whose images are spanned by $x,y,z$, respectively. Then it is easy to see that we have $ [x] \subset  [y] +  [z]$ if and only if $\{ Q,R \}^\perp \subset \{ P  \}^\perp$. Thus, our assumptions on $\phi$ imply that $\xi$ satisfies all the assumptions of the fundamental theorem of projective geometry. Hence, there exists a semilinear bijection $U : H \to H$ such that each rank one projection $P$ is mapped by $\phi$ to the rank one projection whose image is $U ({\rm Im}\, P)$. It is not difficult to see that the condition (\ref{hrvat}) yields that $U$ is actually a unitary or an antiunitary operator multiplied by a nonzero constant. With no loss of generality we can assume that $U$ is a unitary or an antiunitary operator.
Observing that for every $P \in P_1 (H)$ the rank one projection whose  image is $U ({\rm Im}\, P)$ is equal to $UPU^\ast$ we complete the sketch of the proof.

When we speak of applying projective geometry to describe the structure of symmetries we usually mean the reduction of the given problem to Uhlhorn's theorem. Let us illustrate this method by two examples.

\begin{theorem}\label{ludwig}
Let $H$ be a Hilbert space with $\dim H \ge 3$ and $\phi : E(H) \to E (H)$ a bijective map such that (\ref{less}) and (\ref{ortho}) hold true.
Then 
there exists a unitary or antiunitary operator $U : H \to H$ such that $$\phi (A) = UAU^\ast$$
for every $A \in E (H)$.
\end{theorem}

Let us outline the proof of the above statement which is one of the basic results in Ludwig's axiomatic formulation of quantum mechanics. It is clear that the bijectivity of $\phi$ and (\ref{less}) yield $\phi (0) = 0$ and $\phi (I) = I$. Now, if $A$ is a rank one operator in $E(H)$ and $B,C \in E(H)$ satisfy $B,C \le A$ then $B$ and $C$ are comparable. This is a charcateristic property of rank one effects and therefore $\phi$ maps the set of rank one effects onto itself. Since rank one projections are maximal elements in the set of rank one effects we conclude that $\phi (P_1 (H)) = P_1 (H)$. For any pair $P,Q \in P_1 (H)$ we have $P Q = 0$ if and only if $P \le I -Q$ which is equivalent to $\phi (P) \le \phi (I-Q) = I - \phi (Q)$. Thus, $P Q = 0$ holds true if and only $\phi (P) \phi (Q) = 0$. Applying Uhlhorn's theorem we see that there is a unitary or antiunitary operator $U : H \to H$ such that $\phi (P) = UPU^\ast$ for every $P \in P_1 (H)$. After replacing $\phi$ by the map $A \mapsto U^\ast \phi (A) U$ we can assume with no loss of generality that $\phi (P) = P$, $P \in P_1 (H)$.

The reduction to Uhlhorn's theorem shows that $\phi$ has a nice behaviour on the set of all projections of rank one. The rest of the proof is then rather easy. One can see that for every $P \in P_1 (H)$ there exists a bijective increasing function $f_P : [0,1] \to [0,1]$ such that $\phi (tP) = f_P (t)P$, $0 \le t \le 1$. We further have $\phi ((1/2)I) = \phi (I - (1/2)I) = I - \phi ((1/2)I)$, and consequently, $\phi ((1/2)I) = (1/2)I$. It follows from (\ref{less}) that $\phi ((1/2)P) = (1/2)P$ for every $P \in P_1 (H)$. Using (\ref{less}) we can further deduce that $\phi (R) = R$ for every projection $R$ of rank two. The next step is to use the already obtained facts together with the order preserving property to verify that for every pair $P,Q \in P_1 (H)$ satisfying $PQ=0$ we have $\phi ((1/2)P + Q) = (1/2)P + Q$. Now, for every rank one projection $R$ and every real number $s$, $1/2 \le s \le 1$, we can find a pair of projections $P,Q$ such that $tR \le (1/2)P + Q$ if and only if $t \le s$. Thus, $f_R(t) R  \le (1/2)P + Q$ if and only if $t \le s$. This implies that $f_R (s) = s$, $1/2 \le s \le 1$. It is then not difficult to see that we actually have $\phi (tP) = tP$ for every rank one projection $P$ and every real number $t$, $0 \le t \le 1$. It is known that for any two effects $A,B \in E(H)$ we have $A=B$ if and only if 
$$
\max \{ t \, : \, tP \le A \} = \max \{ t \, : \, tP \le B \}
$$
for every rank one projection $P$. It follows that for every $A \in E(H)$ we have $\phi (A) = A$, as desired.

As another example illustrating our method we will outline the proof of the description of symmetries of the set of all projections with respect to commutativity. To make our discussion simpler we will limit ourselves to the finite-dimensional case.
We denote by $P_{n} \subset H_n$ the set of all $n \times n$ projections,
$$
P_{n} = \{ P \in H_n \, : \, P^2 = P = P^\ast \},
$$
and by $P_{n}^1$ the subset of all $n \times n$ projections of rank one. 
Observe first that for every $P \in P_n$ and every $A \in H_n$ the matrices $P$ and $A$ commute if and only if the projection $I-P$ commutes with $A$. It follows that each $\phi : P_n \to P_n$ with the property that for every $P \in P_n$ we have
$$
\phi ( \{ P, I-P \} ) = \{ P, I-P \}
$$
is a symmetry of $P_n$ with respect to commutativity. 

We are now ready to formulate the last result in this section.

\begin{theorem}\label{commutproj}
Let $n \ge 3$ and assume that $\phi : P_n \to P_n$ is a bijective map such that for every pair $P,Q \in P_n$ we have
$$
PQ = QP \iff \phi(P) \phi(Q) = \phi(Q) \phi (P).
$$
Then there exists an $n\times n$ unitary matrix $U$ such that either
$$
\phi ( \{ P, I-P \} ) = \{ UP U^\ast, I-UPU^\ast \}, \ \ \ P \in P_n,
$$
or 
$$
\phi ( \{ P, I-P \} ) = \{ UP^t U^\ast, I-UP^t U^\ast \}, \ \ \ P \in P_n .
$$
\end{theorem}

We outline the proof of this statement that is based on the reduction to Uhlhorn's theorem. We already know that when dealing with such maps $\phi$ there is no distincion between $P$ and $I - P$. Thus, we should have first introduced the equivalence transformation on $P_n$ defined by $P \sim Q$ if $P=Q$ or $P+Q=I$, then defined the commutativity relation on the quotient set in the natural way, and finally defined a new map on the quotient space induced by $\phi$. As we are not interested in the detailed proof but only in main ideas we will ignore this step but we will be aware that if we speak of a rank one projection $P$ we are actually thinking of the pair consisting of the rank one projection $P$ and the projection $I-P$ of corank one.

For any set $\mathcal{P} \subset P_n$ we denote by $\mathcal{P}'$ the commutant of $\mathcal{P}$ in $P_n$, that is, the set of all projections $Q$ that commute with every element of $\mathcal{P}$,
$$\ 
\mathcal{P}' = \{ Q \in P_n \, : \, QP = PQ \ \, {\rm for} \ \, {\rm every} \ \, P \in \mathcal{P} \}.
$$
The commutant of $\mathcal{P}'$ is called the second commutant of $\mathcal{P}$. It is denoted by $\mathcal{P}''$,  $\mathcal{P}''  =  (\mathcal{P}')'$.
Of course, for every subset $\mathcal{P} \subset P_n$ we have
\begin{equation}\label{glasglup}
\phi (\mathcal{P}' ) = \phi(\mathcal{P})' \ \ \ {\rm and} \ \ \  \phi (\mathcal{P}'' ) = \phi(\mathcal{P})'' .
\end{equation}

It is well known and very easy to verify that two hermitian matrices $A$ and $B$ commute if and only if they are simultaneously unitarily similar to diagonal matrices. Thus, if $P$ is a projection of rank one  and $Q$ any projection that commutes with $P$, then up to a unitary similarity we have 
$$
P = \left[ \begin{matrix} 1 & 0 & 0 \cr 0 & 0 & 0 \cr 0 & 0 & 0 \cr \end{matrix} \right] \ \ \ {\rm and} \ \ \ Q = \left[ \begin{matrix} q & 0 & 0 \cr 0 & I & 0 \cr 0 & 0 & 0 \cr \end{matrix} \right],
$$
where $q \in \{ 0,1 \}$. Of course, the second row and the second column, or the third row and the third column in $Q$ may be absent. It is easy to see that
the commutant $\{ P,Q \}'$ consists of all projections of the form
$$
\left[ \begin{matrix} r & 0 & 0 \cr 0 & R_1 & 0 \cr 0 & 0 & R_2 \cr \end{matrix} \right] ,
$$
where $r \in \{ 0,1 \}$ and $R_1$ and $R_2$ are any projections of the appropriate sizes. From here it is easy to conclude the second commutant $\{ P,Q \}''$ consists of all matrices of the form
$$
\left[ \begin{matrix} r_1 & 0 & 0 \cr 0 & r_2 I & 0 \cr 0 & 0 & r_3 I \cr \end{matrix} \right],
$$
where $r_1, r_2 , r_3 \in \{ 0,1 \}$. Thus, the second commutant $\{ P,Q \}''$ has at most 8 elements. 

If, on the other hand, $P$ is a projection that is neither of rank one nor of corank one then we can find a projection $Q$ that commutes with $P$ such that up to unitary similarity we have
$$
P = \left[ \begin{matrix} I & 0 & 0 & 0 \cr 0 & I & 0 & 0 \cr 0 & 0 & 0 & 0 \cr 0 & 0 & 0 & 0 \cr \end{matrix} \right] \ \ \ {\rm and} \ \ \ Q = \left[ \begin{matrix} I & 0 & 0 & 0 \cr 0 & 0 & 0 & 0 \cr 0 & 0 & I & 0 \cr 0 & 0 & 0 & 0 \cr \end{matrix} \right]
$$
and then one can see in the same way as above that the set  $\{ P,Q \}''$ contains 16 elements. 

The trivial projections $0$ and $I$ are the only elements of $P_n$ whose commutant is the whole set $P_n$. It follows that $\phi (\{ 0,I \}) = \{ 0, I \}$. The previous two paragraphs show that we can chacterize rank one projections (that we identify with projections of corank one) among all nontrivial projections with the commutativity relation. Using (\ref{glasglup}) we may assume with no loss of generality that $\phi$ maps the set of projections of rank one onto itself (of course, here we have to take care of our identification of each rank one projection $P$ with $I-P$, but this is just a minor technical difficulty). Thus, restricting $\phi$ to $P_{n}^1$ we arrive at a bijection of $P_{n}^1$ onto itself that preserves commutativity. In the next step we observe that for every pair $P,Q \in P_{n}^1$, $P \not=Q$, we have
$$
PQ = QP \iff PQ =0.
$$
Thus, the restriction of $\phi$ to $P_{n}^1$ satisfies all the assumptions of Uhlhorn's theorem and consequently, after composing $\phi$ with an appropriate unitary slimilarity and the transposition, if necessary, we can assume with no loss of generality that $\phi(P) = P$ for every $P \in P_{n}^1$. Let $Q \in P_n$ be any projection. Then for every $P \in P_{n}^1$ we have
$$
PQ = QP \iff P\phi (Q) = \phi (Q) P,
$$
which is easily seen to imply $\phi (Q) \in \{ Q, I-Q \}$. This completes the sketch of the proof.

\section{Reduction technique}

In the previous two sections we have seen that many problems of describing the general form of given symmetries can be reduced to the description of adjacency preservers or the fundamental theorem of projective geometry.
These two special reduction techniques have been so far the most frequently used general methods when studying symmetries. 

More generally, whenever studying symmetries one of the first questions that we may ask is whether a given problem can be reduced to some other problem of a similar nature. It turns out that many structural results on symmetries have been obtained in this way. Let us illustrate this by three examples.

We will start with the optimal version of Uhlhorn's theorem in the finite-dimensional case. We first assume that $H$ is a separable Hilbert space with $\dim H \ge 3$. Wigner's unitary-antiunitary theorem has been improved in many directions due to its fundamental importance in the mathematical foundations of quantum mechanics. The non-bijective version of Wigner's theorem (for a very simple and elementary proof we refer to \cite{Ge0}) states that if $\phi : P_1 (H) \to P_1 (H)$ is a map (no injectivity or surjectivity is assumed) satisfying (\ref{mcda}), then there exists a linear or conjugate-linear isometry $U : H \to H$ such that $\phi (P) = UPU^\ast$ for every $P \in P_1 (H)$.

We next recall Gleason's theorem, another important result in mathematical physics.  A density operator $D : H \to H$ is defined to be a positive linear operator whose trace is equal to $1$. 
A subset $\{ P_1 , P_2 , \ldots \} \subset P_1 (H)$ is called a complete orthogonal system of projections of rank one, COSP,
if $P_i \perp P_j$ whenever $i \not=j$ and there is no rank one projection $Q$ that is orthogonal to each $P_j$.
Gleason's theorem states that if  $\varphi : P_1 (H) \to [0,1]$ is a function such that for every COSP
$\{ P_1 , P_2 , \ldots \}$ we have
$$
\sum_j \varphi (P_j) = 1,
$$
then there is a density operator   $D : H \to H$ such that
$$
\varphi (P) = {\rm tr}\, (DP)
$$
for every $P \in P_1 (H)$.

Both Wigner's theorem and Gleason's theorem deal with rank one projections (pure states).
While Wigner's theorem is an almost direct consequence of the fundamental theorem of projective geometry, Gleason's theorem is mathematically much deeper and the proof is more involved. It is then natural to ask if Wigner's theorem can be deduced  
from Gleason's theorem. And if this is the case, can we obtain in such a way an improvement of Wigner's theorem?

The answer is positive. We have
the following improvement of Uhlhorn's theorem in the finite-dimensional case \cite{FKKS, PaV}.

\begin{theorem}\label{moremor}
Let $n$ be an integer, $n \ge 3$. Assume that $\phi : P_{n}^1 \to P_{n}^1$ is a map such that for every pair $P,Q \in P_{n}^1$ we have
\begin{equation}\label{topaneki}
PQ = 0 \Rightarrow \phi (P) \phi (Q) = 0.
\end{equation}
Then there exists an $n \times n$ unitary matrix $U$ such that either
$$
\phi (P) = UPU^\ast
$$
for every $P \in P_{n}^1$, or
$$
\phi (P) = UP^t U^\ast
$$
for every $P \in P_{n}^1$.
\end{theorem}

This is an optimal version of Wigner's theorem in the finite-dimensional case. Indeed, like in Uhlhorn's improvement we replace the assumption (\ref{mcda}) by the weaker assumption that only the zero transition probability is preserved. In fact, we go one step further and replace the asumption (\ref{hrvat}) by the weaker assumption (\ref{topaneki}). And we do not assume any regularity condition like injectivity or surjectivity. Let us just mention that in order to get an optimal version of Wigner's theorem in the infinite-dimensional case we need to add a certain rather weak assumption on COSP systems, see \cite{Se3}.

It is clear that under the assumptions of our theorem every COSP is mapped onto some COSP. We take any density matrix $D$ (a positive matrix with trace one) and 
consider the map $\varphi_D : P_{n}^1 \to [0,1]$ defined by
$$
\varphi_D (P) = {\rm tr}\, (D\phi (P)), \ \ \ P \in P_{n}^1 .
$$
By Gleason's theorem we see that for every density matrix $D$ there exists a density matrix $E$ such that
$$
{\rm tr}\, (D\phi (P)) = {\rm tr}\, (EP), \ \ \ P \in P_{n}^1.
$$
In particular, choosing $D = \phi (Q)$ for some fixed $Q \in P_{n}^1$ we get
$$
{\rm tr}\, (\phi (Q)\phi (P)) = {\rm tr}\, (E_Q P), \ \ \ P \in P_{n}^1,
$$
for some density matrix $E_Q$. Set $P=Q$ to conclude that
$$
{\rm tr}\, (E_Q Q) = 1.
$$
It is easy to guess and not difficult to prove that the above equality implies 
$E_Q = Q$.

Hence, for every $P \in P_{n}^1$ we have
$$
{\rm tr}\, (\phi (Q)\phi (P)) = {\rm tr}\, (Q P).
$$
But $Q$ was chosen arbitrarily. Therefore the desired conclusion is a straightforward consequence of
the non-bijective version of Wigner's theorem.

We continue with symmetries of bounded observables with respect to compatibility. Once again we will limit ourselves to the finite-dimensional case. Then the structural result for such symmetries given in the introduction reads as follows.

\begin{theorem}\label{komut}
Let $n$ be an integer, $n \ge 3$, and $\phi : H_n \to H_n$
a bijective map such that for each pair of hermitian matrices $A,B \in H_n$ we have
$$
AB = BA \iff \phi (A) \phi (B) = \phi (B) \phi (A).
$$
Then there exists an $n \times n$ unitary matrix $U$ and for every $A \in H_n$ there is an injective real valued function $f_A$ defined on the spectrum of $A$ such that either
$$
\phi (A) = U f_A (A) U^\ast 
$$
for every $A \in H_n$, or
$$
\phi (A) = U f_A (A^t) U^\ast 
$$
for every $A \in H_n$.
\end{theorem}

We will call the maps appearing in the conclusion of the above theorem standard commutativity preservers.
When proving this theorem it is easy to guess that one can reduce it to Theorem \ref{commutproj}. Let us explain briefly how to do this. For any $A \in H_n$ we denote by $A'$ the commutant of $A$ in $H_n$,
$$
A' = \{ B \in H_n \, : \, AB = BA \}.
$$
It is rather easy to see that for $A,B \in H_n$ we have $A' = B'$ if and only if there exist a unitary matrix $W$ and real numbers $p_1 , \ldots , p_r, q_1 , \ldots, q_r$ such that $p_i \not= p_j$ and $q_i \not= q_j$, $i,j \in \{ 1, \ldots , r\}$, $i \not= j$, and
$$
A = W \,  \left[ \begin{matrix} p_1 I & 0 & \ldots & 0 \cr 0 & p_2 I & \ldots & 0 \cr \vdots & \vdots & \ddots & \vdots \cr 0 & 0 & \ldots & p_r I \cr \end{matrix} \right] \, W^\ast \ \ \ {\rm and} \ \ \
B = W \,  \left[ \begin{matrix} q_1 I & 0 & \ldots & 0 \cr 0 & q_2 I & \ldots & 0 \cr \vdots & \vdots & \ddots & \vdots \cr 0 & 0 & \ldots & q_r I \cr \end{matrix} \right] \, W^\ast .
$$
Here, $I$ and $0$ stand for identity matrices and zero matrices of the appropriate sizes. In this case the commutant $A'$ consists of all hermitian matrices of the form
$$
W \,  \left[ \begin{matrix} \ast & 0 & \ldots & 0 \cr 0 & \ast  & \ldots & 0 \cr \vdots & \vdots & \ddots & \vdots \cr 0 & 0 & \ldots & \ast \cr \end{matrix} \right] \, W^\ast,
$$
where the $\ast$'s stand for any hermitian matrices of the corresponding sizes. And obviously, $B = f(A)$, where $f$ is the injective real valued function defined on the spectrum of $A$ given by $f(p_j) = q_j$, $j= 1, \ldots, r$.

It is now clear that for $A \in H_n$ we have $A' = H_n$ if and only if $A$ is a scalar matrix, that is, $A=tI$ for some real number $t$. Since $\phi (A') = \phi (A)'$, $A \in H_n$, we see that $\phi$ maps the set of scalar matrices onto itself.

The next important observation is that hermitian matrices with two eigenvalues can be characterized as matrices with maximal commutants among all nonscalar matrices. Indeed, assume that $A \in H_n$ has at least three eigenvalues. Then
 there exist a unitary matrix $W$ and real numbers $p_1 , \ldots , p_r$ such that $r \ge 3$, $p_i \not= p_j$, $i,j \in \{ 1, \ldots , r\}$, $i \not= j$, and
$$
A = W \,  \left[ \begin{matrix} p_1 I & 0 & \ldots & 0 \cr 0 & p_2 I & \ldots & 0 \cr \vdots & \vdots & \ddots & \vdots \cr 0 & 0 & \ldots & p_r I \cr \end{matrix} \right] \, W^\ast.
$$ 
But then $B \in H_n$ given by
$$
B = W \,  \left[ \begin{matrix} p_2 I & 0 & \ldots & 0 \cr 0 & p_2 I & \ldots & 0 \cr \vdots & \vdots & \ddots & \vdots \cr 0 & 0 & \ldots & p_r I \cr \end{matrix} \right] \, W^\ast 
$$
is a nonscalar matrix, and $A' \subset B'$, and $A' \not= B'$. It follows that the set of matrices with exactly two eigenvalues is mapped by $\phi$ onto itself.

We introduce an equivalence relation on $H_n$ defined by $A \sim B \iff A' = B'$. It is clear that in each equivalence class consisting of matrices with two eigenvalues there are exactly two projections, $P$ and $I-P$. Thus, the above arguments show that $\phi$ induces a bijective map on the set of all projections (once again we need to deal with the technical detail that $P \sim I-P$ and that $\phi$ does not distinguish these two projections) that preserves commutativity. Hence, we can apply 
Theorem \ref{commutproj} and after composing $\phi$ with an appropriate standard commutativity preserver we can assume with no loss of generality that $\phi (P) = P$ for every $P \in P_n$. Using the fact that for every $A \in H_n$ and every $P \in P_n$ we have
$AP = PA \iff \phi (A) P = P \phi (A)$ we easily conclude that $\phi (A) = f_A (A)$ for some injective real valued  function $f_A$ defined on the spectrum of $A$.

Finally, let $H$ be a Hilbert space and let us consider symmetries of $E(H)$ with respect to coexistency. This problem is quite difficult and we will go through the main steps of the proof of the description of such maps very briefly. For the details we refer to \cite{GeS}. For $A,B \in E(H)$ we write $A \sim B$ if $A$ and $B$ are coexistent. Further, we introduce the notation $A^\sim = \{ C\in E(H)\, : \,  C\sim A \}$. It is not difficult to see that for $A \in E(H)$ we have $A^\sim = E(H)$ if and only if $A$ is a scalar operator. It is far from trivial to verify that for a pair of effects $A,B \in E(H)$ we have $A^\sim = B^\sim$ if and only if  $B\in\{A, A^\perp \} = \{ A, I-A \}$ or both $A$ and $B$ are scalar effects.

The crucial step in studying symmetries of $E(H)$ with respect to coexistency is the introduction of the following somewhat complicated relation on $E(H)$. For $A,B \in E(H)$ we write
$A \prec B$ if and only if for every nonscalar effect $C \in A^\sim$ we can find a nonscalar effect $D \in B^\sim$ such that $C^\sim \subset D^\sim$.
Clearly, $B \prec B$ and $B^\perp \prec B$. And then one can characterize projections in the following way. A nonscalar effect $A \in E(H)$ is a projection if and only if
$${\rm card}\,  \{ B \in  E(H)\setminus \{ tI \, : \, 0 \le t \le 1 \}  \, : \, B \prec A \} = 2.$$
From here we conlude that $\phi$ maps the set of nontrivial projections onto itself. On the set of projections the relation of coexistency coincides with commutativity.
Thus, we have reduced the original problem to Theorem \ref{commutproj}. Therefore we know that $\phi$ has a nice behaviour on the set of all projections. From here one can see that for every symmetry of $E(H)$ with respect to coexistency
there exist a unitary or antiunitary operator $U : H \to H$ and a bijective function $g : [0,1] \to [0,1]$ such that for every nonscalar effect $A$ we have
$$
		\{\phi(A), \phi(A^\perp)\} = \{UAU^\ast, UA^\perp U^\ast\} 
$$
and
$$
		\phi (tI) = g(t)I , \ \ \  t\in [0,1].
	$$


\begin{thebibliography}{99}


\bibitem{BLPY}P. Busch, P. Lahti, J.-P. Pellonp\" a\" a, and K. Ylinen, {\em Quantum measurement}, Springer, 2016.

\bibitem{Dav}E.B. Davies, {\em Quantum theory of open systems}, Academic Press, 1976.




\bibitem{Fa}C.-A. Faure, An elementary proof of the fundamental theorem of projective
geometry, {\em Geom. Dedicata} {\bf 90} (2002), 145--151.

\bibitem{FKKS}A. Fo\v sner, B. Kuzma, T. Kuzma, and N.-S. Sze,
Maps preserving matrix pairs with zero Jordan product,
{\em Linear Multilinear Algebra} {\bf 59} (2011), 507--529. 



\bibitem{Ge0}G.P. Geher, An elementary proof for the non-bijective version of Wigner's theorem, {\em Physics Letters A} {\bf 378} (2014), 2054--2057.


\bibitem{GeS}G. P. Geher and P. \v Semrl, Coexistency on Hilbert space effect algebras and a characterisation of its symmetry transformations, 
{\em Comm. Math. Phys.} {\bf 379} (2020), 1077--1112. 

\bibitem{Hu1}L.K. Hua, Geometries of matrices I. Generalizations of von Staudt's
theorem, {\em Trans. Amer. Math. Soc.} {\bf 57} (1945), 441--481.


\bibitem{Hu2}L.K. Hua, Geometries of matrices ${\rm I}_1$. Arithmetical
construction, {\em Trans. Amer. Math. Soc.} {\bf 57} (1945), 482--490.

\bibitem{Hu3}L.K. Hua, Geometries of symmetric matrices over the real field I,
{\em Dokl. Akad. Nauk. SSSR} {\bf 53} (1946), 95--97.

\bibitem{Hu4}L.K. Hua, Geometries of symmetric matrices over the real field II,
{\em Dokl. Akad. Nauk. SSSR} {\bf 53} (1946), 195--196.

\bibitem{Hu5}L.K. Hua, Geometries of matrices II. Study of involutions in the geometry
of symmetric matrices,
{\em Trans. Amer. Math. Soc.} {\bf 61} (1947), 193--228.

\bibitem{Hu6}L.K. Hua, Geometries of matrices III. Fundamental theorems in the geometries
of symmetric matrices,
{\em Trans. Amer. Math. Soc.} {\bf 61} (1947), 229--255.

\bibitem{Hu7}L.K. Hua, Geometry of symmetric matrices over any field with characteristic
other than two,
{\em Ann. Math.} {\bf 50} (1949), 8--31.

\bibitem{Hu8}L.K. Hua, A theorem on matrices over a sfield and its applications,
{\em Acta Math. Sinica} {\bf 1} (1951), 109--163.

\bibitem{HuS} W.-l. Huang and P. \v Semrl, Adjacency preserving maps on hermitian matrices, {\em Canad. J. Math.} {\bf 60} (2008), 1050--1066.


\bibitem{Kraus}K. Kraus, {\em States, Effects, and Operations},  Lecture Notes in Physics \textbf{190}, Springer-Verlag, 1983.


\bibitem{LudI}G. Ludwig, {\em Foundations of Quantum Mechanics, Vol. I}, (Translated from the German by Carl A. Hein), Springer-Verlag, 1983.

\bibitem{LudII}G. Ludwig, {\em Foundations of Quantum Mechanics, Vol. II}, (Translated from the German by Carl A. Hein), Springer-Verlag, 1985.

\bibitem{Mo1}L. Moln\' ar,  Order-automorphisms of the set of bounded observables,  {\em J. Math. Phys.} {\bf 42} (2001), 5904--5909. 

\bibitem{Mol}L. Moln\' ar, Order automorphisms on positive definite operators and a few applications,
{\em Linear Algebra Appl.} {\bf 434} (2011), 2158 - 2169.

\bibitem{MoK}L. Moln\' ar and E. Kov\' acs, An extension of a characterization of the automorphisms
of Hilbert space effect algebras, {\em Rep. Math. Phys.} {\bf 52} (2003), 141--149.

\bibitem{MoS}L. Moln\' ar and P. \v Semrl, Nonlinear commutativity preserving maps on self-adjoint operators, {\em Quart. J. Math.} {\bf 56} (2005), 589--595.


\bibitem{MS}M. Mori and P. \v Semrl, Loewner's theorem for maps on operator domains, {\em Canad. J. Math.} {\bf 75}  (2023), 912--944.

\bibitem{PaV}M. Pankov and T. Vetterlein, A geometric approach to Wigner-type theorems, {\em Bull. London Math. Soc.} {\bf 53} (2021), 1653--1662.


\bibitem{Se0}P. \v Semrl, Comparability preserving maps on Hilbert space effect algebras, {\em Comm. Math. Phys.} {\bf 313} (2012), 375--384.

\bibitem{Se2}P. \v Semrl, Symmetries on bounded observables - a unified approach based on adjacency
preserving maps,  {\em Integral Equations Operator Theory} {\bf 72} (2012), 7-66.


\bibitem{Se1}P. \v Semrl, Symmetries of Hilbert space effect algebras, {\em J. London Math. Soc.}
{\bf  88}  (2013), 417 - 436.

\bibitem{Se3}P. \v Semrl,
Wigner symmetries and Gleason’s theorem,
{\em J. Phys. A} {\bf 54} (2021), Article ID 315301, 6 p.  






\end{thebibliography}
\end{document}